\documentclass[12pt]{article}

\usepackage{epsf,psfrag,amsmath,amsfonts, amssymb}

\begin{document}
\bibliographystyle{plain}

\newtheorem{theorem}{Theorem}[section]
\newtheorem{lemma}[theorem]{Lemma}
\newtheorem{proposition}[theorem]{Proposition}
\newtheorem{corollary}[theorem]{Corollary}
\newtheorem{conjecture}[theorem]{Conjecture}
\newtheorem{definition}[theorem]{Definition}
\newtheorem{example}[theorem]{Example} 
\newtheorem{condition}{Condition}
\newtheorem{main}{Theorem}
\setlength{\parskip}{\parsep}
\setlength{\parindent}{0pt}

\def \outlineby #1#2#3{\vbox{\hrule\hbox{\vrule\kern #1% 
\vbox{\kern #2 #3\kern #2}\kern #1\vrule}\hrule}}%
\def \endbox {\outlineby{4pt}{4pt}{}}%

\newenvironment{proof}
{\noindent{\bf Proof\ }}{{\hfill \endbox
}\par\vskip2\parsep}

\hfuzz10pt

\newcommand{\cl}[1]{{\mathcal{C}}_{#1}}
\newcommand{\var}{\rm{Var\;}}
\newcommand{\tends}{\rightarrow \infty}
\newcommand{\ep}{{\mathbb {E}}}
\newcommand{\pr}{{\mathbb {P}}}
\newcommand{\re}{{\mathcal{R}}}
\newcommand{\I}{\mathbb {I}}
\newcommand{\rt}{\widetilde{\rho}}
\newcommand{\ro}{\overline{\rho}}
\newcommand{\bds}{\begin{displaystyle}}
\newcommand{\eds}{\end{displaystyle}}
\newcommand{\vc}[1]{{\mathbf{#1}}}
\newcommand{\ra}[2]{{#1}^{(#2)}}
\newcommand{\vd}[1]{\boldsymbol{#1}}

\title{{Fisher Information} inequalities \\
 and the {Central Limit Theorem}}
\author{Oliver Johnson \\ Statslab, Cambridge University \and 
Andrew Barron \\ Dept of Statistics, Yale University}
\date{\today}
\maketitle
\pagestyle{headings}

\begin{abstract} \noindent
We give conditions for an $O(1/n)$ rate of convergence of Fisher  
information and relative entropy in the Central Limit Theorem. We use
the theory of projections in $L^2$ spaces and Poincar\'{e} inequalities, 
to provide a better understanding of the decrease in Fisher Information 
implied by results of Barron and Brown.
We show that if the standardized Fisher Information ever
becomes finite then it converges to zero.
\renewcommand{\thefootnote}{}
\footnote{OTJ is a Fellow of Christ's College, Cambridge, who helped support 
two trips to Yale University during which this paper was written.}
\renewcommand{\thefootnote}{\arabic{footnote}}
\setcounter{footnote}{0}
\end{abstract}

%Proposed Running head: Fisher Information and the CLT
%\renewcommand{\baselinestretch}{1.7}
%\setlength{\baselineskip}{1.4\baselineskip}
%\setlength{\parskip}{1.4\parsep}
%\setlength{\jot}{9pt}
%\setlength{\abovedisplayskip}{11pt}
%\setlength{\belowdisplayskip}{11pt}

\section{Introduction} \label{sec:intro}
Bounds on Shannon entropy and Fisher information have long
been used in proofs of central limit theorems (CLTs), based on 
quantification of the change in information as a result of convolution, 
as in the papers of Linnik (1959)\nocite{linnik}, 
Shimizu (1975)\nocite{shimizu}, Brown (1982)\nocite{brown}, Barron
(1986)\nocite{barron} and Johnson (2000)\nocite{johnson}.  
Each of these papers have a final step 
involving completeness or uniform integrability in which a limit is taken 
without explicitly bounding the information distance from the normal 
distribution.  

The purpose of the present paper is to provide an explicit rate of
convergence of information distances, under certain 
natural conditions on the random variables.  Let $X_1, X_2, ..., X_n$ be 
independent identically distributed random variables with mean $0$, 
variance $\sigma^2$ 
and density function $p(x)$, satisfying Poincare conditions 
(relating $L^2$ norms of mean zero functions to $L^2$ norms of the 
derivative), and let $\phi_{\sigma^2}(x)$ be the corresponding $N(0,\sigma^2)$ 
density.  The relative entropy distance is $D(X)=\int p(x) \log 
\left( p(x)/ \phi_{\sigma^2}(x) \right) dx$.
In the case of random variables with differentiable
densities, the Fisher information distance is 
$J(X)= \sigma^2 \ep[((d/dx)\log p(X)-(d/dx)\log \phi_{\sigma^2}(X))^2]$ 
which is related
to the Fisher information $I(X)=\ep[(d/dx \log p(X))^2]$ via $J(X)=\sigma^2 
I(X)-1$. This is an $L^2$ norm between derivatives of log-densities, and
gives a natural measure of convergence, stronger than existing theorems,
as described in Lemma \ref{lem:compbds}. Note that the quantities
$D$ and $J$ are scale-invariant, that is $D(aX) = D(X)$ and $J(aX) = J(X)$ 
for all non-zero $a$.

Let $U_n=(X_1+ \ldots + X_n)/\sqrt{\sigma^2 n}$ be the standardized sum of the
random variables.  We show that $D(U_n) \le 2R D(U_1)/n\sigma^2$ for all random
variables with Poincare constant $R$, and that 
$J(U_n) \le 2R J(U_1)/n \sigma^2$ 
for all random variables with density function
satisfying a weak differentiability condition. 

The present paper builds on ideas in past work which we
briefly review  here.  In examination of the Fisher information a
central role is played by the score function 
$\rho(y)=(d/dy) \log p(y)=p'(y)/p(y)$.  The score function of the
sum of independent random variables can be expressed in terms of the
score function of the individual random variables, via a conditional
expectation, as has been used in demonstration of convolution inequalities
for Fisher information and Shannon entropy (in the work of Stam 
(1959)\nocite{stam},
Blachman (1965)\nocite{blachman}, and others).  

In particular, if $Y_1$ and  $Y_2$ are
independent and identically distributed with score function $\rho$ then
the score $\ro(u)$ of the sum $Y_1 + Y_2$ is the 
projection of $(\rho(Y_1) + \rho(Y_2))/2$ onto the
linear space of functions of $Y_1 + Y_2$, so by the Pythagorean identity
and rescaling:
\begin{equation}
\frac{I(Y_1) + I(Y_2)}{2} - I \left( \frac{Y_1 + Y_2}{\sqrt{2}} \right)  
=  2\ep \left( \ro \left( Y_1 + Y_2 \right) -
\frac{\rho(Y_1) + \rho(Y_2)}{2} \right)^2
\label{eq:pythag}
\end{equation}
(see Lemma \ref{lem:condit} for details). Hence, since Equation 
(\ref{eq:pythag}) is positive,
%$$
%I \left ( \frac{Y_1 + Y_2}{\sqrt{2}} \right ) \le I(Y_1)
%$$
one deduces that the Fisher Information is decreasing on the powers-of-two
subsequence $S_{k} = U_{2^k}$. Equation (\ref{eq:pythag}) then quantifies the
drop in information $I(S_{k}) - I(S_{k+1})$. 
Identification of the normal as the limiting distribution arises from
examining the difference sequence $I(S_k) - I(S_{k+1})$.

Papers by
Shimizu (1975)\nocite{shimizu}, Brown (1982)\nocite{brown} and
Barron (1986)\nocite{barron} quantify the change in Fisher Information
with each doubling of the sample size, deducing convergence to the normal 
distribution along the powers-of-two subsequence, and convergence of the
whole information sequence $I(U_n)$, by subadditvity of $nI(U_n)$. 
However, these papers only ever consider the behaviour of the 
Fisher Information for $X \sim Y + Z_{\tau}$ (for $Z_\tau$ a small
normal perturbation).

However, in general, we can conclude that if the Fisher Information
$I(S_k)$ is ever finite, since it is
decreasing and bounded below, this difference sequence tends to zero.
Thus the interest is in random variables $Y_1$, $Y_2$ with score
functions for which Equation (\ref{eq:pythag}) is small. This expression
 measures the squared $L^2$ difference between a 
`ridge function' (a function of the sum $Y_1 + Y_2$) and an additive
function (a function of the form $g_1(Y_1) + g_2(Y_2)$). From calculus,
in general, the only functions $f(y_1,y_2) = g_1(y_1) + g_2(y_2)$ that are
both ridge and additive are the linear functions $g_1(y_1) = ay_1 +b$,
$g_2(y_2) = ay_2 +b$, with $a,b_1,b_2$ constants, that is, the functions
for which the derivatives $g_i'(y)$ are constant and equal.

Previous work, as in Lemma 3.1 of 
Brown (1982)\nocite{brown}, (see also Barron (1986)\nocite{barron})
established:
\begin{lemma} \label{lem:brown} For any functions $f$ and $g$ there 
exist some $a$,$b$ such that:
$$ \ep (g(Y_1) - aY_1 -b)^2 \leq \ep \left( f(Y_1 + Y_2) - g(Y_1) - g(Y_2) 
\right)^2, $$
when $Y_1,Y_2$ are independent identically distributed normals.
\end{lemma}
Brown takes $g \in L^2(\phi)$ and considers the
projection $f(s) = \ep( g(Y_1) + g(Y_2) | Y_1 + Y_2 = s)$. For $Y_1, Y_2$
normal, the eigenfunctions of this projection are the Hermite polynomials,
so he can use expansions in this orthogonal Hermite basis.

The main technique used in the present paper will generalize Lemma 
\ref{lem:brown} to a wider class of random variables $Y_1,Y_2$.
For example, consider any $Y_1$ and $Y_2$ IID with finite
Fisher Information $I$. Proposition \ref{prop:main} implies that
given a differentiable ridge function $f(y_1 + y_2)$, with 
closest additive function $g$, then for a certain constant $\mu$:
\begin{equation}
\ep (g'(Y_1) - \mu)^2 
\leq  I \ep \left( f \left( Y_1 + Y_2 \right) - 
g(Y_1) - g(Y_2) \right)^2.  \label{eq:mainineq} \end{equation}
Our (basis-free) proof starts
with $f(Y_1 + Y_2)$, finds its additive part with $g(y_1) =
\ep_{Y_2} f(y_1+Y_2)$ 
and recognises that $g'(y_1) = - \ep_{Y_2} f(y_1 + Y_2) \rho(Y_2)$. 
A Cauchy-Schwarz inequality completes the proof
as detailed in Section \ref{sec:proof}. 

Hence if Equation (\ref{eq:pythag}) is small then $\rho$
is close to a function with derivative close to constant in 
$L^2(Y_1,Y_2)$. 
%
%Equation (\ref{eq:mainineq}) shows that
%if $f(Y_1 + Y_2)$ is close to additive in $L^2(Y_1,Y_2)$
%then $g'$ is close to constant in $L^2$, so $g$ is 
%close to linear as measured by the (Sobolev) $L^2$ norm on
%derivatives. 
Now Poincar\'{e} inequalities 
provide a relationship between $L^2$ norms on functions and the %(Sobolev)
$L^2$ norms on  derivatives:
\begin{definition} \label{def:poinconst}
Given a random variable $Y$, define
the Poincar\'{e} constant $R_Y$:
$$ R_Y = \sup_{g \in H_1(Y)} \frac{ \ep g^2(Y)}{ \ep g'(Y)^2},$$
(where $H_1(Y)$ is the space of absolutely continuous functions $g$ 
such that $\var g(Y) >0$, $\ep g(Y) =0$ and $\ep g^2(Y) < \infty$), and
the restricted Poincar\'{e} constant $R^*_Y$:
$$ R^*_Y = \sup_{g \in H^*_1(Y)} \frac{ \ep g^2(Y)}{ \ep g'(Y)^2},$$
where $H_1^*(Y) = H_1(Y) \cap \{ g: \ep g'(Y)=0 \}$.
\end{definition}
For certain $Y$, $R_Y$ is infinite. However, $R_Y$ is finite
for the normal and other log-concave 
distributions (see for example Klaasen (1985)\nocite{klaasen}, Chernoff
(1981)\nocite{chernoff}, Chen (1982)\nocite{chen2}, 
Cacoullos (1982)\nocite{cacoullos},  
Borovkov and Utev (1984)\nocite{borovkov2}). 
Since we maximise over a smaller set of functions, $R_Y^* \leq R_Y$.
Further,
for $Z \sim N(0,\sigma^2)$, $R^*_Z = \sigma^2/2$, with $g(x) = x^2 -\sigma^2$ 
achieving this 
(we can show this by expanding $g$ in the Hermite basis).

The other important definition that we shall require is that of weak 
differentiability, introduced in Fabian and Hannan (1977)\nocite{fabian}. 
Brown and Gajek (1990)\nocite{brown2} and Lehmann and Casella 
(1998)\nocite{lehmann}) discuss this condition, and provide easier to check
conditions under which it will hold.
\begin{definition}
A random variable $Y$ has weakly differentiable density $p$ if
there exists a function $\rho \in L^2(p)$
such that for all $f$ with $\ep f(Y+u)^2$ finite, the 
function $g(u)=\ep f(Y+u)$ has a derivative $g'(u)$ equal to 
$-\ep[f(Y+u)\rho(Y)]$.
\end{definition}

This is a mild technical condition, allowing an exchange of limit and 
integration. %Indeed $g(u) = \int f(s) p(s-u) ds$
%and its derivative is $g'(u) = -\int f(s) p'(s-u) ds = - \ep f(Y+u) \rho(Y)$.
%
To see the relation to standard differentiability, we can 
take $f(x) = \I( x \in [a,b])$. Then $g(u) = F(b-u) - F(a-u)$ (where 
$F$ is the distribution function of $Y$), and $g'(u) = - \int_{a-u}^{b-u}
p(y) \rho(y)$.
Thus, for any $a,b$ where the distribution function $F$ is differentiable,
$p(b) - p(a) = \int_a^b p(x) \rho(x) dx$.

Using this, one can extend
the Brown inequality Lemma \ref{lem:brown}
to hold (with a constant depending on $I(Y_1)$ and $R_{Y_1}$) 
for a wider class of random variables than just normals.
Since linear score functions correspond to the family of
normal distributions, Equations (\ref{eq:pythag}) and 
(\ref{eq:mainineq}) provide a means to prove the following Central Limit 
Theorems:
\begin{theorem} \label{thm:o1n}
Given $X_1, X_2, \ldots$ IID and with finite variance 
$\sigma^2$, define the normalized sum
$U_n = (\sum_{i=1}^{n} X_i)/\sqrt{n \sigma^2}$.

If $X_i$ are weakly differentiable with 
finite restricted Poincar\'{e} constant $R^*$ then
$$ J(U_n) \leq \frac{2R^*}{n \sigma^2}J(X) \mbox{ for all $n$}.$$
If $X_i$ have finite Poincar\'{e} constant $R$, then
$$ D(U_n) \leq \frac{2R}{n \sigma^2} D(X) \mbox{ for all $n$}.$$
\end{theorem}
\begin{proof}
See Sections \ref{sec:proof} and \ref{sec:respoin} for the proof of the Fisher
information bound.
Notice that for $X$ normal, $2R^* = \sigma^2$, 
so the `closer to normal $X$ is', the closer the bound becomes to $J(X)/n$. 
 
The relative entropy bound is a corollary.
Using an integral form of the de Bruijn identity
(Lemma 1 of Barron (1986)\nocite{barron}), the relative entropy 
$D(X)$ can be
expressed as an integral of $J(\sqrt{1-t} X + \sqrt{t}Z)$ (that is, a linear
combination of $X$ and a standard normal $Z$). Now, if $X$ has finite
Poincar\'{e} constant $R$, then for each $t$, the $(\sqrt{1-t} X + \sqrt{t}Z)$
itself has Poincar\'{e} constant $\leq R$, so by Theorem \ref{thm:o1n},
the result follows. \end{proof}

This $O(1/n)$ rate of convergence is perhaps to be expected. For example if
$X_i$ is exponentially distributed, and hence $U_n$ has a
$\Gamma(n)$ distribution, then $J(U_n) = 2/(n-2)$, consistent with this.
In fact, by extending the Cram\'{e}r-Rao lower bound we deduce that
\begin{lemma} 
If $\ep X^4$ is finite, then 
$$ \liminf_{n \rightarrow \infty} \; nJ(U_n) \geq s^2/3,$$
where $s$ is the skewness, $m_3(X)/m_2(X)^{3/2}$ (writing $m_r(X)$ 
for the centred $r$th moment of $X$).
\end{lemma}
\begin{proof}
For any function f, the positivity of $\ep(\rho_U(U) +f(U))^2$ implies that
$\ep \rho_U(U)^2 \geq \ep(2f'(U) - f(U)^2)$, giving a whole family of bounds. 
Assume that $\ep U =0$
and take $ f(u) = (u - u^2 m_3(U)/m_4(U))/m_2(U)$.
%Then $\ep f'(U) = 1/m_2$, and $\ep f(U)^2 = a^2 m_4 + 2a m_3/m_2 + 1/m_2$.
%The optimal value of $a$ is $-m_3/(m_2 m_4)$.
% 
%Hence 
Then $J(U) = m_2(U) \ep \rho_U(U)^2 -1 
\geq m_3(U)^2/(m_2(U) m_4(U))$. 
Now, since $m_2(U_n) = m_2(X)$, $m_3(U_n) = m_3(X)/\sqrt{n}$ and
$m_4(U_n) = m_4(X)/n + 3m_2(X)^2(n-1)/n$, the result follows.
\end{proof} 

Further, this $O(1/n)$ convergence is consistent with 
estimates of Berry--Esseen type
which give a $O(1/\sqrt{n})$ rate of weak convergence. 
The following lemma shows the relationship between 
convergence in Fisher Information, and several weaker forms of convergence:
\begin{lemma} \label{lem:compbds}
If $X$ is a random variable with density $f$, and $\phi$
is a standard normal, then:
\begin{eqnarray*}
\sup_x |f(x) - \phi(x)| & \leq & 
\left( 1 + \sqrt{\frac{6}{\pi}} \right) \sqrt{J(X)}, \\
\int |f(x) - \phi(x)| dx \leq
2d_H(f,\phi) & \leq & \sqrt{2} \sqrt{J(X)}, 
\end{eqnarray*}
where $d_H( f ,\phi)$ is the Hellinger 
distance $\left( \int | \sqrt{f(x)} - \sqrt{\phi(x)}|^2 dx \right)^{1/2}$.
\end{lemma} 
\begin{proof} The first bound comes from Shimizu (1975)\nocite{shimizu}. The 
second inequality tightens a bound of Shimizu. Since:
$$ \sqrt{\phi(x)} \frac{\partial}{\partial x} \sqrt{\frac{f(x)}{\phi(x)}} 
= \frac{1}{2} \left( \frac{f'(x)}{\sqrt{f(x)}} + x \sqrt{f(x)} \right),$$
we deduce from the Poincar\'{e} inequality for $\phi$ that:
$$ J(X) = 4 \int \phi(x) 
\left( \frac{\partial}{\partial x} \sqrt{\frac{f(x)}{\phi(x)}} \right)^2
\geq 4 \int \phi(x) 
\left( \sqrt{\frac{f(x)}{\phi(x)}} - \mu \right)^2 
= 4(1 - \mu^2),$$
where $\mu = \ep_\phi \sqrt{f/\phi}$, so $d^2_H(f,\phi) = 2(2 - 2\mu)
\leq 4(1-\mu^2)$.
\end{proof}

Recent work by Ball et al (2002)\nocite{ball} has also considered the rate of 
convergence of these quantities. Their paper obtains similar results,
but by a very different method, involving transportation costs and a 
variational characterisation of Fisher Information.

Unfortunately, Poincar\'{e} constants are not finite for all 
distributions $Y$. 
Indeed, as Borovkov and Utev (1984)\nocite{borovkov2} point out, 
if $R_Y < \infty$, then by considering $g_n(x) = |x|^n$, we inductively 
deduce that all the moments of $Y$ are finite. 
From the Berry-Esseen Theorem we know
that only $(2+\delta)$th moment conditions are enough to ensure an explicit 
$O(1/n^{\delta/2})$ rate of weak convergence, for $0 < \delta \leq 1$. 
In Section \ref{sec:nopoin} we describe a proof of
Fisher Information convergence under only second moment conditions,
though without an explicit rate. This is an extension of Barron's
Lemma 2, which only holds for random variables with a normal perturbation.
\begin{theorem} \label{thm:fishconv}
Given $X_1, X_2, \ldots$ weakly differentiable IID with finite variance 
$\sigma^2$, define the normalized sum
$U_n = (\sum_{i=1}^{n} X_i)/\sqrt{n \sigma^2}$.
If $J(U_m)$ is finite for some $m$ then $$ \lim_{n \tends} J(U_n) = 0.$$
\end{theorem}
{\bf Note:} This extends Lemma 2 of Barron (1986)\nocite{barron}, which only 
holds when $X$ is of the form $Y + Z_\tau$.
\section{Projection of functions in $L^2$} \label{sec:proof}
Although the main application of the following Proposition
will concern score functions,
we present it as an abstract result concerning projection of functions
in $L^2(Y_1,Y_2)$.
\begin{proposition} \label{prop:main}
Consider independent random variables $Y_1, Y_2$ with weakly 
differentiable densities and functions $f, h_1, h_2$ such that 
$ \ep[(f(Y_1+Y_2))^2]$ is finite
and 
$\ep f(Y_1  + Y_2) = 0$. We find functions $g_1,g_2$ and a 
constant $\mu$ such that for any $\beta \in [0,1]$:
\begin{eqnarray*} 
\lefteqn{ \ep \left( f(Y_1 + Y_2) - h_1(Y_1) - h_2(Y_2) \right)^2 } \\
& \geq & \ep (g_1(Y_1) - h_1(Y_1))^2 + \ep (g_2(Y_2) - h_2(Y_2))^2 \\
& & + (\overline{I})^{-1} 
\left( \beta \ep \left( g_1'(Y_1) - \mu \right)^2
+ (1-\beta) \ep \left( g_2'(Y_2) - \mu \right)^2 \right),
\end{eqnarray*}
where $\overline{I} = (1-\beta) I(Y_1) + \beta I(Y_2)$.
\end{proposition}
\begin{proof}
Firstly, given $
\ep \left( f(Y_1 + Y_2) - h_1(Y_1) - h_2(Y_2) \right)^2$, we can
replace $h_1$ and $h_2$ by functions $g_1, g_2$ which reduce it even further.
This follows since this $L^2$ distance is minimised over choices of $h_1$,
$h_2$ by considering the projections, which remove the additive part:
\begin{eqnarray*}
g_1(u) & = & \ep_{Y_2}  f( u + Y_2),  \\
g_2(v) & = & \ep_{Y_1}  f( Y_1 + v), 
\end{eqnarray*}
so the Pythagorean relation tells us that the LHS equals
$$ \ep (g_1(Y_1) - h_1(Y_1))^2 + \ep (g_2(Y_2) - 
h_2(Y_2))^2 + \ep \left( f(Y_1 + Y_2) - g_1(Y_1) - g_2(Y_2) \right)^2.$$
Having removed the additive part of $f$, we hope that what remains will be
small in magnitude. Hence, the inner product of what remains 
and certain functions of the variables should be small.  Specifically we 
define
\begin{eqnarray*}
r_1(u) & = & \ep_{Y_2} \left[ \left( f( u + Y_2) - g_1(u) - g_2(Y_2) \right) 
\rho_2(Y_2) \right], \\
r_2(v) & = & \ep_{Y_1} \left[ \left( f( Y_1 + v) - g_1(Y_1) - g_2(v) \right) 
\rho_1(Y_1) \right],
\end{eqnarray*}
and show that we can control their norms.
Indeed, by Cauchy-Schwarz, for any $u$:
$$ r^2_1(u) \leq  \ep_{Y_2} 
\left( f( u + Y_2) - g_1(u) - g_2(Y_2) \right)^2 \ep \rho^2_2(Y_2),$$
so taking expectations over $Y_1$, we deduce that
\begin{equation} \label{eq:f1bound}
 \ep r^2_1(Y_1) \leq  \ep 
\left( f( Y_1 + Y_2) - g_1(Y_1) - g_2(Y_2) \right)^2 I(Y_2).
\end{equation}
Similarly,
\begin{equation} \label{eq:f2bound}
\ep r^2_2(Y_2) \leq  \ep
\left( f( Y_1 + Y_2)- g_1(Y_1) - g_2(Y_2) \right)^2 I(Y_1).
\end{equation}
The assumption that $\ep f(Y_1+Y_2)^2$ is finite implies that for almost
  every $u$ we have $\ep f(u+Y_2)^2$ finite.  Consider any such $u$.
  The weak differentiability of $p_2$ (the density of $Y_2$) gives that the
  function $g_1(u)=
 \ep[f(u+Y_2)]$ has derivative $g_1^{'}(u) = -\ep f(u+Y_2)
\rho_2 (Y_2)$.
  Also weak differentiability trivially yields $\ep \rho_2 (Y_2)=0$, so setting
  $\mu= - \ep g_2(Y_2) \rho_2 (Y_2)$, we recognize that $r_1(u)$ defined above
  simplifies to
$$ r_1(u) = - \left( g_1'(u) - \mu \right)$$
Using the similar expression for  
$r_2(v) = - (g_2'(v) - \mu)$, and adding 
$\beta$ times Equation (\ref{eq:f1bound}) to $(1-\beta)$ times 
Equation (\ref{eq:f2bound}), we deduce the result.
\end{proof}

{\bf Note:} this inequality holds in general, for any weakly
differentiable $Y_1,Y_2$ with finite
Fisher Information, whereas previous such expressions have only held in the
case of $Y_i \sim U_i + Z_\tau$, for some $U_i$.

{\bf Note:} this inequality allows for independent
random variables that are not identically distributed. Armed with it, one may
provide Central Limit Theorems giving information convergence to the 
normal for random variables satisfying a uniform Lindeberg-type condition
(see also Johnson (2000)\nocite{johnson}). In certain cases we can provide
a rate of convergence. 

{\bf Note:}
we can produce a similar expression using a similar method for 
finite-dimensional random
vectors $\vc{Y}_1,\vc{Y}_2$. Weak differentiability can be defined in this
case, and $\rho_i = (\partial/\partial x_i)(\log p(x))$ 
will usually be 
the $i$th component of the score vector function $\vd{\rho}$.
Similar analysis in this case can lead to an alternative proof of the 
Theorems in Johnson and Suhov (2001)\nocite{johnson3}.
\section{Rate of convergence} \label{sec:respoin}
If $Y_1, Y_2$ have finite restricted Poincar\'{e} constants $R^*_1, 
R^*_2$ then we can extend 
Lemma \ref{lem:brown} from the case of normal $Y_1,Y_2$ to more
general distributions, providing an explicit
exponential rate of convergence of Fisher Information.
We can apply Proposition \ref{prop:main} 
because the score functions of sums can be expressed
as $L^2$ projections.
\begin{lemma} \label{lem:condit}
Let $S=Y_1+Y_2$ where $Y_1$ and $Y_2$ are independent and $Y_2$ is
  weakly differentiable with score function $\rho_2$.  Then $S$ is weakly
  differentiable with score function 
$$ \ro(s) = \ep[\rho_2(Y_2)|S=s].$$
Hence for independent weakly differentiable
random variables $Y_1$ and $Y_2$, with score 
functions $\rho_1$ and $\rho_2$, writing $\ro$ for the score function of $S$
$$
\frac{I(Y_1) + I(Y_2)}{2} - I \left(
\frac{Y_1 + Y_2}{\sqrt{2}} \right)  
=  2 \ep \left( \ro \left( S \right) -
\frac{\rho_1(Y_1) + \rho_2(Y_2)}{2} \right)^2.$$
\end{lemma}
\begin{proof}
For any square integrable test function $T(u+S)$,
define $T_2 (v)= \ep[T(v+Y_1)]$ so that $T_2 (u+Y_2) = \ep[T(u+S)|Y_2]$.
Note that $\ep[(T_2(u+Y_2))^2] \le \ep[T^2 (u+S)] < \infty$, 
so that weak differentability can be applied to it. That is
  we have $\ep [T(u+S) \ro(S)]= \ep [T(u+S) \rho_2(Y_2)]=
\ep [T_2(u+Y_2) \rho_2(Y_2)]
  =-(d/du) \ep T_2(u+Y_2) = -(d/du)\ep T(u+S)$.

Then if both random variables are weakly differentiable,
$\ro = \ep[(\rho_1(Y_1) + \rho_2(Y_2))/2 |S=s]$.
Thus by the Pythagorean identity, the result follows, since we can 
rescale: $\rho_{aX}(x) = \rho_X(x/a)/a$ and $J(aX) = J(X)/a^2$. \end{proof}

\begin{proposition} \label{prop:fishrate}
Consider $Y_1,Y_2$ IID and weakly differentiable with variance $\sigma^2$ and
restricted Poincar\'{e} constant $R^*$. Then
$$ J \left( \frac{Y_1 + Y_2}{\sqrt{2}} \right) \leq J(Y_1) 
\left( \frac{2R^*}{\sigma^2+2R^*} \right).$$
\end{proposition}
\begin{proof}
Without loss of generality, suppose $Y_i$ have mean 0 and variance 1,
since we can just rescale, using $R_{aX}^* = a^2 R_{X}^*$.
Write $J$ and $I$ for the standardardized and non-standardized Fisher
Information of $Y$, and $J'$ and $I'$ for the corresponding quantities
for $(Y_1+Y_2)/\sqrt{2}$.
By rescaling Lemma \ref{lem:condit}, for projections $g$, writing 
$\rt$ for the score function of $(Y_1 + Y_2)/\sqrt{2}$: 
\begin{eqnarray*}
\lefteqn{J(Y_1) - J \left(
\frac{Y_1 + Y_2}{\sqrt{2}} \right)}  \\
& = & \ep \left( \rt \left( \frac{Y_1 + Y_2}{\sqrt{2}} \right) - \frac{
g(Y_1) + g(Y_2)}{\sqrt{2}} \right)^2 + \ep( \rho_1(Y_1) - g(Y_1))^2.
\end{eqnarray*}
\begin{figure}[ht!]
\begin{center}
\begin{psfrags}
\psfrag{a}[][]{$-(Y_1 + Y_2)/\sqrt{2}$}
\psfrag{b}[l][]{$\rt\left((Y_1 + Y_2)/\sqrt{2} \right)$}
\psfrag{c}[l][]{$(g(Y_1) + g(Y_2))/\sqrt{2}$}
\psfrag{d}[l][]{$(\rho(Y_1) + \rho(Y_2))/\sqrt{2}$}
\psfrag{e}[l][]{$(h(Y_1) + h(Y_2))/\sqrt{2}$}
\epsfysize=8.9cm
\leavevmode\epsfbox{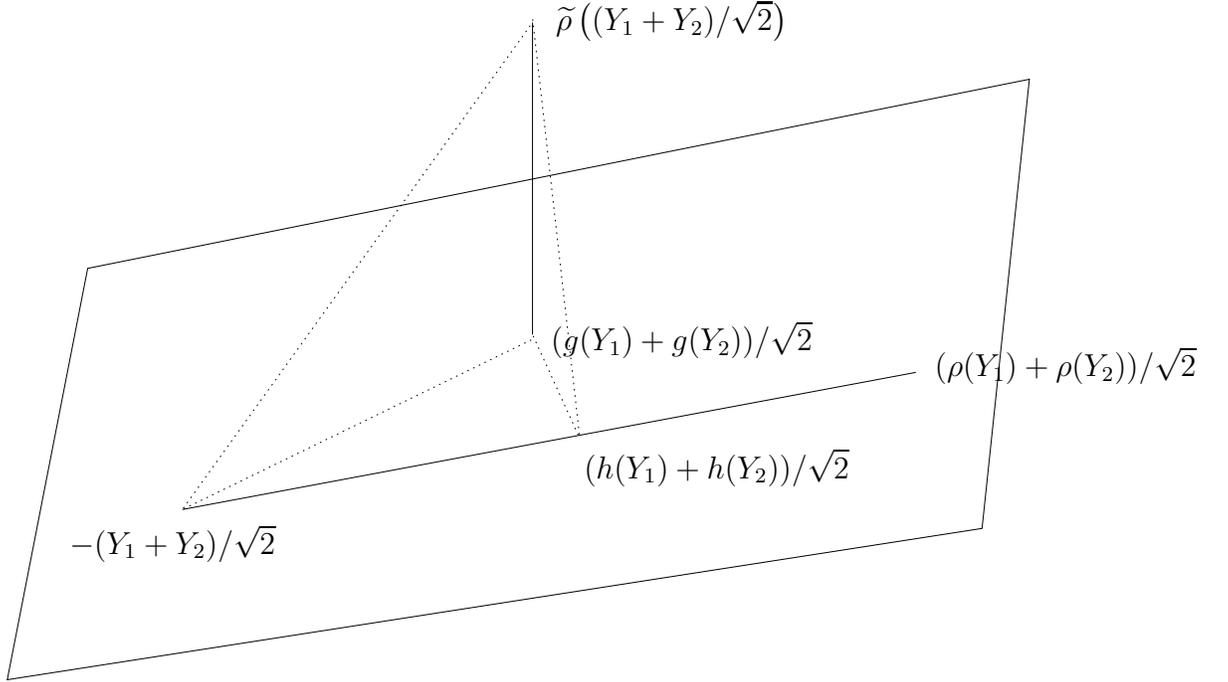}
\caption{Role of projections}
\label{pic:proj}
\end{psfrags}
\end{center}
\end{figure}
Now, consider the projection of $\rt$ into the space of additive 
functions, shown as a plane
in Figure \ref{pic:proj}, where $(h(Y_1) + h(Y_2))/\sqrt{2}$ 
is the closest point to $\rt$ on the line between $-(Y_1 + Y_2)/\sqrt{2}$ 
and $(\rho(Y_1) + \rho(Y_2))/\sqrt{2}$, so 
that $\ep (g(Y) +Y)^2 \geq \ep (h(Y) + Y)^2$.

Further, we know that $h$ corresponds to the value of $\lambda$ which 
minimises:
$$ \ep \left( \rt \left( \frac{Y_1 + Y_2}{\sqrt{2}} \right)
-  \left( \lambda \left( \frac{\rho(Y_1) + \rho(Y_2)}{\sqrt{2}} \right) 
+ (1-\lambda) \left( \frac{Y_1 + Y_2}{\sqrt{2}} \right) \right) \right)^2.$$
Since in general $\ep (U - \lambda V)^2$ is minimised at 
$\lambda = \ep UV/\ep V^2$,
in this case the minimising $\lambda = J'/J$, so $h$ is $J'/J$ of the way
along the line. This tells us that 
$\ep (h(Y) + Y)^2 = (J'/J)^2 \ep (\rho(Y) + Y)^2 = J^{'2}/J$. 

Overall then, we deduce that $\ep (g(Y) +Y)^2 \geq J^{'2}/J$, and by
Pythagoras,
$\ep \left( \rt \left( (Y_1 + Y_2)/\sqrt{2} \right) - (g(Y_1) + g(Y_2))
/\sqrt{2} \right)^2 \leq J' - J^{'2}/J$.

Now applying Proposition \ref{prop:main} to the first bracket, we can see
that the factor of $I$ in the denominator that Proposition \ref{prop:main}
implies will actually cancel, simplifying the expression.
\begin{eqnarray*}
J' - J^{'2}/J 
& \geq & \ep \left( \rt \left( \frac{Y_1 + Y_2}{\sqrt{2}} \right) - \frac{
g(Y_1) + g(Y_2)}{\sqrt{2}} \right)^2 \\
& \geq & \frac{\ep (g_1'(Y_1) - \mu)^2}{2I} 
\geq  \frac{\ep (g_1(Y_1) - \mu Y_1)^2}{2R^*I} \\
& = &  \frac{\ep (g_1(Y_1) + Y_1)^2 + (-\mu-1)^2}{2R^*I} 
\geq   \left( \frac{J^{'2}}{J}  + J^{'2} \right) \frac{1}{2R^*I}
=   \frac{J^{'2}}{2R^*J}
\end{eqnarray*}
since $\ep g_1'(Y_1) - \mu = 0$ where $\mu = -I'$,
and since $\ep Y_1 g_1(Y_1) = - 1$, 
so rearranging, we obtain the result.

Note that $g_1$ will be absolutely continuous, so we can apply the Poincar\'{e}
inequality. This follows since 
$ \int_v^w g_1'(u) du = - \int_v^w \int f(s) p'(s-u)ds du$
and Fubini's Theorem tells us that is $- \int f(s) \int_v^w p'(s-u) du ds
= - \int f(s) (p(s-w) - p(s-v)) ds = g_1(w) - g_1(v)$,
as the random variables have a density everywhere.
\end{proof}
A more careful analysis generalises Proposition \ref{prop:fishrate},
to obtain Theorem \ref{thm:o1n} by performing 
successive projections onto smaller additive spaces.  
For a given function $f$, 
define a series of functions by $ f_n = f$, and for $m <n$:
$$ f_m \left( \frac{ X_1 + \ldots + X_m}{\sqrt{n}} \right) =
\ep_{X_{m+1}} f_{m+1} \left( \frac{ X_1 + \ldots + X_m + X_{m+1}}{\sqrt{n}} 
\right).$$
Further, define $g(u) = \sqrt{n} \ep f
\left( \frac{ X_1 + \ldots + X_{n-1} + u}{\sqrt{n}}  \right)$.
At step $i$, we approximate the function $f$ by 
$f_i ((X_1+...+X_i)/\sqrt{n})$
plus a sum of $g(X_j)$ for $j>i$, which is the best approximation
onto the linear space of such partially additive functions.
\begin{lemma} \label{lem:bd1}
Defining the squared distance between successive projections to be
$$ t_i = \ep \left( 
f_i \left( \frac{ X_1 + \ldots + X_i}{\sqrt{n}} \right)
- f_{i-1} \left( \frac{ X_1 + \ldots + X_{i-1}}{\sqrt{n}} \right) - 
\frac{1}{\sqrt{n}} g(X_i) \right)^2,$$
then for $X_i$ IID and weakly differentiable:
$$ t_i \geq \frac{(i-1)}{nI(X)} \ep (g'(X) - \mu)^2.$$
\end{lemma}
\begin{proof} 
Evaluate the function 
\begin{eqnarray*} 
r(z) & = & \ep \bigg[ \left( 
f_i \left( \frac{ X_1 + \ldots + X_{i-1} + z}{\sqrt{n}} \right)
- f_{i-1} \left( \frac{ X_1 + \ldots + X_{i-1}}{\sqrt{n}} \right) - 
\frac{1}{\sqrt{n}} g(z) \right)  \\
& & \; \; \; \; \; \times  (\rho(X_1) + \ldots \rho(X_{i-1})) \bigg],
\end{eqnarray*}
in two different ways. Firstly, again using weak differentiability of $X_i$:
$$ r(z) = - \left( \frac{i-1}{\sqrt{n}} \right) \left( g'(z) - \mu \right),$$
where $\mu = \ep f_{i-1}' = \ep f'$.
Secondly, we apply Cauchy-Schwarz to $r(z)^2$, and take expected values,
to deduce that
$ \ep r(X)^2 
\leq t_i(i-1)I(X)$. Putting these together, the result follows.
\end{proof}
\begin{lemma} \label{lem:bd2} For $X_i$ IID,
the sum of these squared distances $t_i$ is $s_n = \sum_{i=1}^n t_i$, where
$$s_m = \ep \left( f_m \left( \frac{ X_1 + \ldots + X_m}{\sqrt{n}} \right)
- \sum_{i=1}^m \frac{ g(X_i)}{\sqrt{n}} \right)^2.$$
\end{lemma}
\begin{proof}
Since $s_m = \ep f_m^2 - (m/n) \ep g^2$, and since
$ t_m = \ep f_m^2 - \ep f_{m-1}^2 - (1/n) \ep g^2$, we can rearrange to 
obtain:
$$ s_m = (t_m + \ep f_{m-1}^2 + (1/n) \ep g^2) - (m/n) \ep g^2 
= t_m + s_{m-1},$$
so summing the telescoping sum, the result follows.
\end{proof}
Combining Lemma \ref{lem:bd1} and Lemma \ref{lem:bd2}, we deduce that:
\begin{equation} \label{eq:ct}
 s_n \geq \sum_{i=1}^n \frac{(i-1)}{nI(X)} \ep (g'(X) - \mu)^2
= \frac{(n-1)}{2I(X)} \ep (g'(X) - \mu)^2.\end{equation}
\begin{proof}{\bf of Theorem \ref{thm:o1n}} Again, assuming that $X$ has 
variance 1, and
writing $J'$ for $J(U_n)$, and $J$ for $J(X)$,
as before we know that $\ep (g(X) + X)^2 \geq J^{'2}/J$ and
$s_n = \left( \rho_n - \sum g(X_i)/\sqrt{n} \right)^2 \leq 
J'(1- J'/J)$. Hence by Equation (\ref{eq:ct}), we deduce that:
\begin{eqnarray*}
J'(1-J'/J) & \geq &  s_n \geq \frac{(n-1)}{2I(X)} \ep (g'(X) - \mu)^2 \\
& \geq & \frac{(n-1)}{2R^*I(X)} \ep (g(X) - \mu X)^2 \\
& \geq &  \frac{(n-1)}{2R^*I(X)} \left( \frac{J^{'2}}{J} + J^{'2} \right) 
=   \frac{(n-1)}{2R^*} \frac{J^{'2}}{J}.
\end{eqnarray*}
Thus, in general, rescaling gives: 
$$J(U_n) \leq \frac{2R^*}{2R^* + (n-1)\sigma^2}J(X),$$
and the result follows.\end{proof}
{\bf Note:} for $X$  discrete-valued, 
$X+Z_{\tau}$ has a finite Poincar\'{e} constant, and 
hence this calculation of an explicit rate of convergence 
of $J(S_n + Z_{\tau})$ still holds. Via Lemma \ref{lem:compbds}
we know that $S_n+Z_\tau$ converges weakly for any $\tau$
and hence $S_n$ converges weakly to the standard normal. 
\section{Convergence of Fisher Information} \label{sec:nopoin}
We can still obtain convergence of the Fisher Information, though without
such an attractive rate of convergence, if the Poincar\'{e} constants are not
finite. We will
need uniform control over the tails of the Fisher Information, and then
will bound it on the rest of the region using the projection
arguments of Section \ref{sec:proof}. Recall that for 
$I(X)$ finite, the density of $X$ is bounded
(since $p(y) \leq \int p(x) |\rho(x)| dx \leq \sqrt{I(X)}$).
\begin{definition}
Given a function $\psi$, we define the following class:
$$ {\cal C}_\psi = \{ X: \ep X = 0, \sigma^2 = \ep X^2 < \infty,
\sigma^2 \ep \rho(X)^2 \I(|X| \geq \sigma R) \leq \psi(R) 
\mbox{ for all $R$.}\} $$
\end{definition}
\begin{lemma} \label{lem:tailbd}
For $X_1, X_2, \ldots $ IID with finite variance and finite
$I(X)$, then $U_m \in {\cal C}_{\psi}$ for all $m$
where $\psi(R) = \ep \rho( X)^2  \I(|X| \geq R) + C/R^{1/2}$. \end{lemma}
\begin{proof} We take the common variance to be equal to 1 and
use the notation that $p$ and $\rho$ stand for the density and score
function of a single $X$, and $p_r$ for the density of $X_1 + \ldots X_r$.
We know that $U_m$ has score function 
$ \rho_m(u) = \ep \left( \left. \sum_i \rho(X_i ) \right| U_m = u \right)/
\sqrt{m}$,
so by the conditional version of Jensen's
inequality
\begin{equation} \label{eq:condjen}
\rho_m(u)^2 \leq \ep
( \rho(X_1)^2 | U_m = u)
+ (m-1) \ep (\rho(X_1) \rho(X_2) | U_m = u). \end{equation}
Consider the two terms of Equation (\ref{eq:condjen}) separately, 
firstly writing $W$ for $X_2 + \ldots X_m$:
\begin{eqnarray*}
\lefteqn{\ep_{U_m} \ep[ \rho( X_1)^2 | U_m] \I(|U_m| \geq R) } \\
& \leq & \ep \rho( X_1)^2 \left( \I(|X_1| \geq R, |U_m| \geq R)
+ \I(|X_1| < R, |U_m| \geq R) \right) \\
& \leq & \ep \rho( X_1)^2 \left( \I(|X_1| \geq R)
+ \I(|W| \geq R(\sqrt{m}-1) ) \right) \\
& \leq & \ep \rho( X)^2  \I(|X| \geq R)
+ \frac{I(X) (m-1)}{R^2(\sqrt{m}-1)^2}
\end{eqnarray*}
Then for any $u$:
\begin{eqnarray*}
\lefteqn{\ep (\rho(X_1) \rho(X_2) | U_m = u)} \\
& = & \iint  \frac{\sqrt{m} p(v) p(w) 
p_{m-2}(u\sqrt{m} - v - w)}{q_m(u)} \rho(v) \rho(w) dv dw \\
& = & \frac{\sqrt{m}}{q_m(u)} 
\int p_{m-2}(u\sqrt{m} - x) 
\int \frac{\partial p}{\partial v}(v) \frac{\partial p}{\partial x}(x-v) dv dx
  \\
%& = & \frac{\sqrt{m}}{q_m(u)} \int p_{m-2}(u\sqrt{m} - x) p''_2(x) dx  \\
%& = & \frac{\sqrt{m}}{q_m(u)} p''_m(u\sqrt{m})   
%=  \frac{q''_m(u)}{mq_m(u)}.   
\end{eqnarray*}
So the second term of Equation (\ref{eq:condjen}) is $q_m'(-R) - q_m'(R)$ and
we need a function $\psi'$ such that for all $R$:
\begin{equation} |q'_m(R)| \leq \psi'(|R|) \label{eq:derbd} \end{equation}
For all $m$, 
$q_m(x) \leq \sqrt{I(U_m)} \leq \sqrt{I}$ for some $I$.
Since for any $m$:
\begin{eqnarray*}
q_{2m}'(u) & = & 2 \int q_m'(v) q_m( u \sqrt{2} - v) dv \\
& \leq & 2^{3/4} \left( \int \frac{q_m'(v)^2}{q_m(v)} dv \right)^{1/2}
\left( \int \sqrt{2} q_m(v) q_m^2( u \sqrt{2} - v) dv \right)^{1/2} \\
& \leq & 2^{3/4} \sqrt{I} 
\left( \sqrt{I} \int \sqrt{2} q_m(v) q_m( u \sqrt{2} - v) dv \right)^{1/2} \\
& \leq & (2I)^{3/4}\sqrt{q_{2m}(u)}
\end{eqnarray*}
(a similar bound will hold for $q_{2m+1}$) and
$$ q_m(u) \leq  \int_u^{\infty} |q_m'(v)| dv 
\leq  \left( \int_u^{\infty} \frac{q_m'(v)^2}{q_m(v)} dv \right)^{1/2}
\left( \int_u^{\infty} q_m(v) dv \right)^{1/2} \\
 \leq  \frac{\sqrt{I}}{u},$$
we deduce that 
Equation (\ref{eq:derbd}) holds, with $\psi'(R) = 2^{3/4} I/R^{1/2}.$
Note that under a $(2+\delta)$th moment condition, we obtain $\psi'(R) =
C/R^{(2+\delta)/4}$.
\end{proof}
By results of Brown (1982)\nocite{brown}, we know that under a finite variance
condition,  there exists $\theta(R)$  such that  $\ep X^2  \I(|X| \geq
\sigma R)  \leq \theta(R)$.  If  in addition, $\ep  |X|^{2+\delta}$ is
finite for  some $\delta$, the Rosenthal inequality  implies that $\ep
|U_n|^{2+\delta}$ is  uniformly bounded, so  we can take  $\theta(R) =
1/R^{\delta}$.

The other ingredient we require is a bound on the Poincar\'{e} constant 
$R_{U_n}^T$ (the Poincar\'{e} constant of $U_n$ conditioned on $|U_n| \leq T$).
\begin{lemma} \label{lem:poinbd}
If $I(X) $ is finite then there exist $R(T)$ and $N(T)$ such that
for all $T$, $R_{U_n}^T \leq R(T)$ for $n \geq N(T)$.
\end{lemma}
\begin{proof} 
Writing $d_n = \sup_A | f_n(A)
- \phi(A) |$ (which tends to zero), since $f_n$ is bounded then:
\begin{eqnarray*}
\left| f_{2n}(x) - \phi(x) \right|
& \leq & \sqrt{2} \left| \int f_n(\sqrt{2} x -y) (f_n(y) - \phi(y)) dy \right|
\\
& & + \sqrt{2} \left| \int \phi(\sqrt{2} x -y) (f_n(y) - \phi(y)) dy  \right|
 \\
& \leq & \sqrt{2} \left( \| f_n \|_{\infty} + \| \phi \|_{\infty} \right)
\int (f_n(y) - \phi(y)) I(f_n(y) \geq \phi(y)) dy \\
& \leq & \left(\sqrt{2I} + \sqrt{1/\pi} \right) d_n
\end{eqnarray*}
Now, for given $T$, take $N(T) = 2 \min \left\{ n : \left(\sqrt{2I} + 
\sqrt{1/\pi} 
\right) d_n \leq \phi(T)/2 \right\}$. This implies that
$f_{n}(x) \geq \phi(T)/2$, 
for $x \in [-T,T]$ and $n \geq N(T)$, so $R(T) = 2/\phi(T)$ means
\begin{eqnarray*}
\int_{-T}^{x} y f_{n}(y) dy & \leq & R(T) f_{n}(x), \mbox{ for 
$0 \geq x \geq -T$} \\
\int_{x}^{T} y f_{n}(y) dy & \leq & R(T) f_{n}(x), \mbox{ for $0 \leq x \leq 
T$},
\end{eqnarray*}
since the LHS is always less than $1$, so by Theorem 1 of Borovkov and Utev 
(1984)\nocite{borovkov2} we are done.
\end{proof}
Combining these two results gives:

\begin{proof}{\bf of Theorem \ref{thm:fishconv}}
Firstly, using the projection inequalities (Proposition \ref{prop:main}),
there exist a function $g$ and constants $\mu, \nu$ such that:
\begin{eqnarray*}
\lefteqn{J(U_n) - J(U_{2n})} \\
& = & \ep \left( \rho_{2n}(U_{2n}) - \frac{1}{\sqrt{2}} \rho_n(U_n)
- \frac{1}{\sqrt{2}} \rho_n(U_n') \right)^2 \\
& \geq & \ep( \rho_n(U_n) - g(U_n))^2  + \frac{1}{2I(U_n)} 
\ep( g'(U_n) - \mu)^2  \\
& \geq & \ep( \rho_n(U_n) - g(U_n))^2 \I(|U_n| \leq T)  
+ \frac{1}{2I(U_n)} \ep( g'(U_n) - \mu)^2 \I(|U_n| \leq T) \\ 
& \geq & \ep( \rho_n(U_n) - g(U_n))^2 \I(|U_n| \leq T)  \\
& &
+ \frac{1}{2R^T_{U_n} 
I(U_n)} \ep( g(U_n) - \mu U_n - \nu)^2 \I(|U_n| \leq T) \\
& \geq & 
\frac{1}{1+2R^T_{U_n} I(U_n)} 
\ep( \rho_n(U_n) - \mu U_n - \nu)^2 \I(|U_n| \leq T) \\
\end{eqnarray*}
Now $\mu = \ep \rho'_{2n}(U_n) = -I(U_{2n})$, and 
$\nu = - \ep (\rho_n + U_n) \I(|U_n| \geq T)$, so
$$ J(U_n) \leq (1+2R^T_{U_n} I(U_n))(J(U_n) - J(U_{2n})) 
+ \ep( \rho_n(U_n) - \mu U_n - \nu)^2 \I(|U_n| > T),$$
and hence by Lemmas \ref{lem:tailbd} and \ref{lem:poinbd}, 
for some function $\zeta(T)$ such that $\zeta(T) \rightarrow 0$ as
$T \tends$:
$$ J(U_n) \leq (1+2R^T_{U_n} I(U_n))
(J(U_n) - J(U_{2n})) + \zeta(T).$$
For any $\epsilon$ we can find $T_0$ such that $\zeta(T_0) \leq \epsilon$, 
for all $n \geq N(T_0)$, then $(1+2R^T_{U_n} I(U_n))(J(U_n) - J(U_{2n}) 
\leq (1+2R(T_0) I)(J(U_n) - J(U_{2n})  \leq \epsilon$ for $n$ sufficiently 
large. \end{proof}

\makeatletter

{\bf \large Addresses} 

O.T.Johnson, Statslab, 
CMS, Wilberforce Road, Cambridge, CB3 0WB, UK. Email: {\tt otj1000@cam.ac.uk}. 

A.R.Barron, Dept of Statistics, Yale University,
PO Box 208290, New Haven, Connecticut 06520-8290, USA.
Email: {\tt andrew.barron@yale.edu}.
\makeatother


\begin{thebibliography}{10}

\bibitem{ball}
Ball, K., Barthe, F. and Naor, A. (2002).
\newblock Entropy jumps in the presence of a spectral gap
\newblock {\em Preprint}.

\bibitem{barron}
Barron, A.R. (1986).
\newblock Entropy and the {Central Limit Theorem}.
\newblock {\em Ann. Probab.} {\bf 14}, 336--342.

\bibitem{blachman} 
Blachman, N.M. (1965)
\newblock The convolution inequality for entropy powers.
\newblock {\em IEEE Trans. Inform. Theory} {\bf 11}, 267--271.

\bibitem{borovkov2}
Borovkov, A.A. and Utev, S.A. (1984).
\newblock On an inequality and a related characterisation of the normal
  distribution.
\newblock {\em Theory Probab. Appl.}{\bf 28}, 219--228.

\bibitem{brown}
Brown, L.D. (1982).
\newblock A proof of the {Central Limit Theorem} motivated by the
  {Cram\'{e}r-Rao} inequality.
\newblock In G.~Kallianpur, P.R. Krishnaiah, and J.K. Ghosh, editors, {\em
  Statistics and Probability: Essays in Honour of {C.R. Rao}}, pages 141--148.
  North-Holland, New York.

\bibitem{brown2}
Brown, L.D. and Gajek, L. (1990)
\newblock Information inequalities for the {B}ayes risk.
\newblock {\em Ann. Statist.}{\bf 18}, 1578--1594.

\bibitem{cacoullos}
Cacoullos, Th. (1982).
\newblock On upper and lower bounds for the variance of a function of a random
  variable.
\newblock {\em Ann. Probab.} {\bf 10}, 799--809.

\bibitem{chen2}
Chen, L.H.Y. (1982)
\newblock An inequality for the normal distribution.
\newblock {\em J. Multivariate Anal.} {\bf 12}, 306--315.

\bibitem{chernoff}
Chernoff, H. (1981)
\newblock A note on an inequality involving the normal distribution.
\newblock {\em Ann. Probab.} {\bf 9}, 533--535 .

\bibitem{fabian}
Fabian, V. and Hannan, J. (1977)
\newblock On the {C}ram\'er-{R}ao inequality.
\newblock {\em Ann. Statist.} {\bf 5}, 197--205.

\bibitem{johnson}
Johnson, O.T. (2000)
\newblock Entropy inequalities and the {Central Limit Theorem}.
\newblock {\em Stochastic Process Appl.} {\bf 88}, 291--304.

\bibitem{johnson3}
Johnson, O.T. and Suhov, Y.M. (2001)
\newblock Entropy and random vectors.
\newblock {\em J. Statist Phys.} {\bf 104}, 147--167.

\bibitem{klaasen}
Klaasen, C.A.J. (1985).
\newblock On an inequality of {Chernoff}.
\newblock {\em Ann. Probab.} {\bf 13}, 966--974.

\bibitem{lehmann}
Lehmann, E. and Casella, G. (1998).
\newblock {\em Theory of point estimation}.
\newblock Springer Texts in Statistics. Springer-Verlag, New York, second
  edition.

\bibitem{linnik}
Linnik, Y.V. (1959)
\newblock An information-theoretic proof of the {Central Limit Theorem} with
  the {Lindeberg Condition}.
\newblock {\em Theory Probab. Appl.} {\bf 4}, 288--299.


\bibitem{shimizu}
Shimizu, R. (1975)
\newblock On {Fisher's } amount of information for location family.
\newblock In G.P.Patil et~al, editor, {\em Statistical Distributions in
  Scientific Work, Volume 3}, pages 305--312. Reidel.

\bibitem{stam} 
Stam, A.J. (1959)
\newblock Some inequalities satisfied by the quantities of information of
  {Fisher} and {Shannon}.
\newblock {\em Inform. and Control} {\bf 2}, 101--112.



\end{thebibliography}
\end{document}